\documentclass{amsart}

\usepackage{amsmath,amssymb,amscd,amsfonts}

\newtheorem{theorem}{Theorem}
\newcommand{\bt}{\begin{theorem}}
\newcommand{\et}{\end{theorem}}
\newtheorem{lemma}{Lemma}
\newcommand{\bl}{\begin{lemma}}
\newcommand{\el}{\end{lemma}}
\newtheorem{corollary}{Corollary}
\newcommand{\bc}{\begin{corollary}}
\newcommand{\ec}{\end{corollary}}
\newcommand{\beq}{\begin{equation}}
\newcommand{\eeq}{\end{equation}}
\newcommand{\benum}{\begin{enumerate}}
\newcommand{\eenum}{\end{enumerate}}

\newcommand{\Z}{\ensuremath{\mathbf Z}}

\newcommand{\mbF}{\ensuremath{ \mathbf F}}

\newcommand{\mca}{\ensuremath{ \mathcal A}}

\newcommand{\mcp}{\ensuremath{ \mathcal P}}

\newcommand{\F}{\ensuremath{\mathbf F }}
\newcommand{\Fq}{\ensuremath{{\mathbf F}_q}}
\newcommand{\mbt}{\ensuremath{\mathbf t}}

\DeclareMathOperator{\card}{\text{card}}

\newcommand{\bmat}{\left(\begin{matrix}}
\newcommand{\emat}{\end{matrix}\right)}

\DeclareMathOperator{\qqand}{\qquad\text{and}\qquad}

\title{The Bose-Chowla argument for Sidon sets}

\author{Melvyn B. Nathanson}
\address{Lehman College (CUNY), Bronx, New York 10468}
\email{melvyn.nathanson@lehman.cuny.edu}



\thanks{Supported in part by a grant from the PSC-CUNY Research Award Program.}

\date{\today}

\begin{document}
\maketitle 

\begin{abstract} 
Let $h \geq 2$ and let ${ \mathcal A} = (A_1,\ldots, A_h)$ be an $h$-tuple of sets of integers.  
For  nonzero integers $c_1,\ldots, c_h$, consider the linear form 
$\varphi =  c_1 x_1 + c_2x_2 + \cdots + c_h x_h$.
The \emph{representation function} $R_{ \mathcal{A},\varphi}(n)$ 
 counts the number of $h$-tuples $(a_1,\ldots, a_h) \in A_1 \times \cdots \times A_h$ such that 
$\varphi(a_1,\ldots, a_h) = n$. 
The $h$-tuple $\mathcal{A}$  is a \emph{$\varphi$-Sidon system of multiplicity $g$} 
if  $R_{\mathcal A,\varphi}(n) \leq g$ for all  $n \in \Z$.  
For every positive integer $g$, let $F_{\varphi,g}(n)$ denote the largest integer $q$ such that there exists 
a $\varphi$-Sidon system $\mathcal {A} = (A_1,\ldots, A_h)$ of multiplicity $g$ with 
\[
A_i \subseteq [1,n] \qqand |A_i| = q
\]
for all $i =1,\ldots, h$.   It is proved that, for all linear forms $\varphi$, 
\[
\limsup_{n\rightarrow \infty} \frac{F_{\varphi,g}(n)}{n^{1/h}}  < \infty 
\]
and, for linear forms $\varphi$ whose coefficients $c_i$ satisfy a certain divisibility condition, 
\[
\liminf_{n\rightarrow\infty} \frac{F_{\varphi,h!}(n)}{n^{1/h}} \geq 1.
\]
\end{abstract}

\section{Classical Sidon sets} 
Let $A$ be a subset of an additive abelian group or semigroup $\Lambda$.  
For every positive integer $h$,  let $A^h$ be the set 
of all $h$-tuples of elements of $A$.  
The symmetric group $S_h$ acts on the set $A^h$ by permutation of coordinates:  
For all $\sigma \in S_h$ and $(a_1,\ldots, a_h) \in A^h$, we define 
\beq                         \label{BC-Sidon:action} 
\sigma \left( a_1,\ldots, a_h \right) = \left( a_{\sigma(1)}, \ldots, a_{\sigma(h)} \right).
\eeq
The orbits of this action define an  equivalence relation $\sim$ on $A^h$:  
For $(a_1,\ldots, a_h) \in A^h$ and $(a'_1,\ldots, a'_h) \in A^h $, we have  
\[
(a_1,\ldots, a_h) \sim (a'_1,\ldots, a'_h)
\]
if and only if there is a permutation $\sigma \in S_h$ such that 
\[
a'_i = a_{\sigma(i)} \qquad \text{for all $i \in \{1,\ldots, h\}$.}
\]
Let $[ a_1,\ldots, a_h ]$ denote the orbit of the $h$-tuple $(a_1,\ldots, a_h)$  
and let 
\[
A^h/S_h = \left\{ [a_1,\ldots, a_h] : (a_1,\ldots, a_h) \in A^h \right\}. 
\]
be the set of equivalence classes of $A^h$.  
The number of $h$-tuples in the orbit  $[a_1,\ldots, a_h ]$ is at most $h!$  
and is equal to $h!$ if and only if the coordinates of the $h$-tuple are distinct.    

Consider the linear form
\[
\varphi  = \varphi (x_1, \ldots, x_h) = x_1 + \cdots + x_h.
\]
If $(a_1,\ldots, a_h) \sim (a'_1,\ldots, a'_h)$, then for some permutation $\sigma \in S_h$ we have 
\[
\varphi(a_1,\ldots, a_h) = a_1+\cdots + a_h = a_{\sigma(1)} + \cdots + a_{\sigma(h)} = \varphi(a'_1,\ldots, a'_h).  
\]
Thus, the function  $[\varphi]$ defined by  
\[
[\varphi]( [a_1,\ldots, a_h] ) =  \varphi(a_1,\ldots, a_h) = a_1+\cdots + a_h 
\]
 is a well-defined function on the orbit space $A^h /S_h$.  

The \emph{$h$-fold sumset} of $A$ is the set  $hA$ of all sums of $h$ not necessarily distinct elements of $A$: 
\begin{align*}
hA & = \left\{ a_1+\cdots + a_h: (a_1,\ldots, a_h) \in A^h  \right\}  \\ 
& = \left\{ \varphi(a_1,\ldots, a_h):  (a_1,\ldots, a_h) \in A^h \right\} \\
& = \left\{ [\varphi] ( [a_1,\ldots, a_h] ):  [a_1,\ldots, a_h] \in A^h/S_h \right\} 
\end{align*}
 
We define two representation functions for the $h$-fold sumset $hA$:   For $w \in \Lambda$, 
\begin{align*}
r_{A,h}(w) & = \card  \left\{ [a_1,\ldots, a_h] \in A^h/S_h : a_1 + \cdots + a_h = w   \right\}  
\end{align*} 
and 
\begin{align*}
R_{A,h}(w) & = \card\left\{ (a_1,\ldots, a_h) \in A^h: a_1 + \cdots + a_h = w \right\}. 
\end{align*} 
Because the coordinates of  the $h$-tuple $(a_1,\ldots, a_h)$ are not necessarily distinct, we have  
\[
R_{A,h}(w) \leq h!  \ r_{A,h}(w)
\]
for all $w \in \Lambda$.

The set $A$ is a \emph{classical Sidon set of order $h$} or a \emph{$B_h$-set} 
if $r_{A,h}(w) \leq 1$ for all $w \in \Lambda$.
The set $A$ is a \emph{classical Sidon set of order $h$ and multiplicity $g$ }
 or, simply, a  \emph{$B_h[g]$-set}  if $r_{A,h}(w) \leq g$ for all $w \in \Lambda$.
 
 For every integer $m \geq 2$ we also define the \emph{modular representation function} 
\[
r^{(m)}_{A,h}(w)  = \card  \left\{ [a_1,\ldots, a_h] \in A^h/S_h : a_1 + \cdots + a_h \equiv w \pmod{m}   \right\}.  
\]
The set $A$ is a \emph{classical Sidon set of order $h$ modulo $m$} 
if $r^{(m)}_{A,h}(w) \leq 1$ for all $w \in \Lambda$.

If $A$ is a classical Sidon set of order $h$ modulo $m$ for some $m \geq 2$, then 
 $A$ is a classical Sidon set of order $h$. 

Halberstam and Roth~\cite{halb-roth66}  
and O'Bryant~\cite{obry04} are excellent surveys of results on Sidon sets obtained before 2004.

\section{The Bose-Chowla argument}
Consider classical Sidon sets for the group \Z\ of  integers. 
Let $F_h(n)$ denote the largest Sidon set of order $h$ contained in the 
set of consecutive integers $\{1,2,\ldots, n\}$.  
A simple counting argument shows that $F_h(n) \ll n^{1/h}$ and so 
\[
\limsup_{n\rightarrow\infty} \frac{F_h(n)}{n^{1/h}} < \infty. 
\]
Bose and Chowla~\cite{bose-chow62}  proved that 
\[
\liminf_{n\rightarrow\infty} \frac{F_h(n)}{n^{1/h}}  > 0.
\]
Their main result is the following.

\bt[Bose-Chowla]         \label{BC-Sidon:theorem:BC-1}
Let $h \geq 2$ and let $q$ be a prime power.   
Let $\Fq = \{  \lambda_1,\ldots, \lambda_{q}\}$ be a finite field with $q$ elements 
and let $\F_{q^h}$ be an extension field of \Fq\ of degree $h$.   
Let $\theta$ be a generator of the cyclic group $\F_{q^h}^{\times}$.  
For all $\lambda_j \in \Fq$, there is a unique  integer $a_j \in \{1,2,\ldots, q^h - 2\}$ such that 
\[
\theta^{a_j} = \theta - \lambda_j.
\]
The set 
\[
A = \left\{ a_j : j \in \{1,\ldots, q\} \right\} \subseteq \{1,2,\ldots, q^h - 2\}
\]
is a classical Sidon set of order $h$ modulo $q^h -1$ and  cardinality $q$. 
\et

\begin{proof}
The element $\theta$ generates the field extention $\mbF_{q^h}/\Fq$, 
and so the minimal polynomial of $\theta$ has degree $h$.  
Because $\theta$ has order $^{q^h-1}$ in the cyclic group $\mbF_{q^h}^{\times}$,  for all integers $u$ and $v$  
we have  $\theta^u = \theta^v$ if and only if $u\equiv v \pmod{q^h-1}$. 

Note that $\theta \notin \Fq$ because $h \geq 2$ and $\theta$ generate $\F_{q^h}^{\times}$.  
Thus,  $\theta - \lambda_j \neq 0$ for all $\lambda_j \in \Fq$.  
It follows that for all $j \in \{1,\ldots, q\}$ there is a unique  integer 
$a_j \in \{0,1,2,\ldots, q^h - 2\}$ such that 
\[
\theta^{a_j} = \theta - \lambda_j.
\]   
If $a_j=0$, then $1 = \theta^0= \theta - \lambda_j$ and 
$\theta = \lambda_j +1 \in \Fq$, which is absurd.  
Therefore, $a_j \in \{1,2,\ldots, q^h - 2\}$.  

We have $a_j = a_k$ 
if and only if $\theta-\lambda_j = \theta^{a_j} = \theta^{a_k} = \theta-\lambda_k$ 
if and only if $j = k$.  It follows that the set 
\[
A = \left\{ a_j : j \in \{1,\ldots, q \}  \right\}
\]
has cardinality $q$. 

Let $(a_1,\ldots, a_h)$ and $(a'_1,\ldots, a'_h)$ be $h$-tuples of elements of $A$ such that 
\[
a_1 \leq a_2 \leq \cdots \leq a_h \qqand a'_1 \leq a'_2 \leq \cdots \leq a'_h. 
\]
There exist unique elements $\lambda_j \in \Fq$ and $\lambda'_j \in \Fq$ such that 
\[
\theta^{a_j} = \theta - \lambda_j \qqand \theta^{a'_j} = \theta - \lambda'_j 
\] 
for all $j \in \{1,2,\ldots, h\}$.  If 
\[
a_1 + a_2 + \cdots + a_h \equiv  a'_1 + a'_2 + \cdots + a'_h \pmod{q^h-1}
\]
then 
\begin{align*} 
\prod_{j =1}^h (\theta - \lambda_j)  & =  \prod_{j =1}^h \theta^{a_j} = \theta^{a_1 + a_2 + \cdots + a_h} 
 = \theta^{a'_1 + a'_2 + \cdots + a'_h}  = \prod_{j=1}^h \theta^{a'_j} \\
 &  = \prod_{j =1}^h (\theta - \lambda'_j).
\end{align*} 

The polynomial 
\[
f(t) = \prod_{j =1}^h (t - \lambda_j) - \prod_{j =1}^h (t - \lambda'_j)
\]
is either the zero polynomial or a nonzero polynomial of degree at most $h-1$ with coefficients in \Fq.  
Moreover, 
\[
f(\theta) = 0. 
\]
The minimal polynomial of $\theta$ has degree $h$, and so $\theta$ is not a root of a nonzero polynomial of degree less than $h$.  
Therefore,  $f(t)$ must be the zero polynomial, and so the polynomials  
$\prod_{j =1}^h (t - \lambda_j)$ and $\prod_{j =1}^h (t - \lambda'_j)$ have the same roots with the same multiplicities.  
Thus, $(\lambda_1,\ldots, \lambda_h)$ is a permutation of  $(\lambda'_1,\ldots, \lambda'_h)$,  and so 
$(a_1,\ldots, a_h) = (a'_1,\ldots, a'_h)$.  
Equivalently, if $(a_1,\ldots, a_h) \in A^h$ and
$(a'_1,\ldots, a'_h) \in A^h$ and $(a_1,\ldots, a_h) \neq (a'_1,\ldots, a'_h)$, then 
\[ 
a_1 + a_2 + \cdots + a_h \not\equiv a'_1 + a'_2 + \cdots + a'_h \pmod{q^h -1}
\]
and so $A$ is a  classical Sidon set of order $h$ modulo $q^h-1$.   
This completes the proof.  
\end{proof}

\bc               \label{BC-Sidon:corollary:BC-1}
For every prime power $q$ and every integer $h \geq 2$, 
\[
F_h(q^h -2 ) \geq q.
\]
\ec

\begin{proof}
A  classical Sidon set of order $h$ modulo $q^h-1$ is a  classical Sidon set of order $h$.   
\end{proof}

\bt[Bose-Chowla]         \label{BC-Sidon:theorem:BC-2}
For  every integer $h \geq 2$, 
\[
\liminf_{n\rightarrow \infty} \frac{F_h(n )}{n^{1/h}} \geq 1.
\]
\et

\begin{proof}
It has been known since Hoheisel~\cite{hohe30} (and, more recently, Heath-Brown~\cite{heat88b}) 
that there is a real number $\alpha$ 
with $0 < \alpha < 1$ such that 
if $p$ and $p'$ are consecutive primes, then $p' -p < p^{\alpha}$.

For every integer $n \geq 2^h$, let $p$ be the largest prime such that
$p \leq n^{1/h}$ and let $p'$ be the smallest prime such that $p' > n^{1/h}$.  
The primes $p$ and $p'$ are consecutive, and so, by Hoheisel, 
\[
p \leq n^{1/h} < p' \leq p+p^{\alpha} \leq p+  n^{\alpha/h} 
\]
The function $F_h(n)$ is increasing.  
Applying Corollary~\ref{BC-Sidon:corollary:BC-1} with $q = p$, we obtain 
\[
F_h(n)  \geq F_h(p^h)  \geq F_h(p^h-2) \geq p 
\geq n^{1/h} - n^{\alpha /h}.
\]
Therefore,
\[
\liminf_{n\rightarrow \infty} \frac{F_h(n )}{n^{1/h}} 
\geq \liminf_{n\rightarrow \infty} \frac{ n^{1/h} - n^{\alpha /h}}{n^{1/h}} 
= \liminf_{n\rightarrow \infty} \left( 1 - \frac{ 1}{ n^{(1-\alpha) /h}} \right) = 1.
\]
This completes the proof. 
\end{proof}

\section{Sidon systems for linear forms}

Fix an integer $h \geq 2$.  
Let $\mca = (A_1,\ldots, A_h)$ be an $h$-tuple of sets of integers and let 
\[
A_1 \times \cdots \times A_h  = \left\{ (a_1,\ldots, a_h) : a_i \in A_i \text{ for } i = 1,\ldots, h \right\}.
\]
For  nonzero integers $c_1,\ldots, c_h$, we consider the linear form 
\[
\varphi =  c_1 x_1 + c_2x_2 + \cdots + c_h x_h
\]
and the set 
\[
\varphi(\mca) 
= \left\{   \varphi(a_1,\ldots, a_h) : (a_1,\ldots, a_h) \in A_1 \times \cdots \times A_h  \right\}.
\]

For every integer $n$, the \emph{representation function} $R_{\mca,\varphi}(n)$ 
 counts the number of $h$-tuples $(a_1,\ldots, a_h) \in A_1 \times \cdots \times A_h$ such that 
$\varphi(a_1,\ldots, a_h) = n$. 

Let $W$ be a finite or infinite set of integers.   
The $h$-tuple \mca\ is a \emph{$\varphi$-Sidon system for $W$} 
if $R_{\mca,\varphi}(w) \leq 1$ for all  $w \in W$.  
If \mca\ is a $\varphi$-Sidon system for $W$, 
then the statements 
\[
(a_1,\ldots, a_h) \in A_1\times \cdots \times A_h, \qquad   (a'_1,\ldots, a'_h) \in A_1\times \cdots \times A_h
\]
 and 
\[
 \varphi(a_1,\ldots, a_h) = \varphi(a'_1,\ldots, a'_h) = w \in W 
\]
imply  
\[
(a_1,\ldots, a_h) =  (a'_1,\ldots, a'_h). 
\]
More generally, the $h$-tuple \mca\ is a \emph{$\varphi$-Sidon system of multiplicity $g$ for $W$} 
if $R_{\mca,\varphi}(w) \leq g$ for all  $w \in W$.  

The $h$-tuple \mca\ is a \emph{$\varphi$-Sidon system} if \mca\ is a \emph{$\varphi$-Sidon system} for \Z, 
that is, if  $R_{\mca,\varphi}(n) \leq 1$ for all integers $n$.  
The $h$-tuple \mca\ is a \emph{$\varphi$-Sidon system of multiplicity $g$} 
if  $R_{\mca,\varphi}(n) \leq g$ for all integers $n$.

Let   $\mca = (A_1,\ldots, A_h)$ be an $h$-tuple of finite sets of integers.  We have  
\[
 \prod_{i=1}^h |A_i| = \sum_{w\in \varphi(\mca)} R_{\mca,\varphi}(w)
\]
If \mca\ is a $\varphi$-Sidon system, then 
\[
 \sum_{w\in \varphi(\mca)} R_{\mca,\varphi}(w) = |\varphi(\mca)|. 
\]
If \mca\ is a $\varphi$-Sidon system of mutiplicity $g$, 
then 
\[
 \sum_{w\in \varphi(\mca)} R_{\mca,\varphi}(w) \leq g|\varphi(\mca)|. 
\]

For example, let $d_1,\ldots,  d_h$ be integers such that $d_i \geq 2$ for all $i = 1,\ldots, h$ 
and let $d = d_1d_2\cdots d_h$.  Let $W = [0, d-1]$.  
Consider the finite sets
\[
A_i = [0, d_i-1] \qquad \text{for $i = 1,\ldots h$ } 
\]
and the linear form 
\[
\varphi = x_1 + d_1x_2 + d_1d_2 x_3 + \cdots + d_1d_2\cdots d_{h-1} x_h.  
\]
We have $R_{\mca,\varphi}(w) = 1$ for all  $w \in W$ 
and  $R_{\mca,\varphi}(w) = 0$ for all integers $w \not\in W$.  
The $h$-tuple $\mca = (A_1, \ldots, A_h)$ is a Sidon system.


\section{The size of $\varphi$-Sidon systems} 

Let $\mca = (A_1,\ldots, A_h)$ be an $h$-tuple of sets of integers 
and let $(t_1,\ldots, t_h)$ be an $h$-tuple of integers.   
For $i \in \{1,\ldots, h\}$, the translate of the set $A_i$ by the integer $t_i$ is the set 
\[
A_i + t_i = \{a_i+t: a_i \in A_i\}
\]
and the translate of  \mca\ by  $(t_1,\ldots, t_h)$ is the $h$-tuple 
\[
\mca + \mbt = \mca+ (t_1,\ldots, t_h)  = (A_1+t_1,\ldots, A_h + t_h).  
\]

Let $\varphi = \sum_{i=1}^h c_ix_i$, where the coefficients $c_i$ are nonzero integers and let 
\[
t^* = \varphi(t_1,\ldots, t_h).
\]

For all $(a_1,\ldots, a_h) \in A_1 \times \cdots \times A_h$, we have 
\begin{align*}
\varphi(a_1+t_1,\ldots, a_h+t_h) 
& = \sum_{i=1}^h c_i (a_i+t_i) = \sum_{i=1}^h c_i a_i + \sum_{i=1}^h c_i t_i \\
& = \varphi(a_1,\ldots, a_h) + \varphi(t_1,\ldots, t_h) \\
 & = \varphi(a_1,\ldots, a_h) + t^*.
\end{align*}
Thus, $\varphi(a_1,\ldots, a_h) = b$ if and only if $ \varphi(a_1+t_1,\ldots, a_h+t_h) = b + t^*$, and so 
\[
R_{\mca,\varphi}(b) = R_{\mca+ \mbt,\varphi}(b+t^*).
\]
It follows that for $h$-tuples of nonempty  finite sets of integers or $h$-tuples  
of nonempty sets of integers that are bounded below, 
it suffices to consider only $h$-tuples of sets of nonnegative integers $(A_1,\ldots, A_h)$ 
with $0 \in A_i$ for all $i \in \{1,\ldots, h\}$.

Let $\varphi = c_1 x_1 + \cdots + c_h x_h$, where the coefficients $c_i$ are nonzero integers, and let 
\[
C = \sum_{i=1}^h |c_i|.
\]
Let $\mca = (A_1,\ldots, A_h)$ be an $h$-tuple of finite sets of integers 
and let $n$ be a positive integer such that $|a_i| \leq n$ for all $a_i \in \bigcup_{i=1}^h A_i$.   
For all $(a_1, \ldots, a_h) \in  A_1 \times \cdots \times A_h$ we have 
\[
\left| \varphi(a_1, \ldots, a_h) \right| = \left| \sum_{i=1}^h c_ia_i\right| 
\leq  \sum_{i=1}^h  \left|c_i\right|  |a_i| \leq C n
\]
and so 
\[
\varphi(A_1,\ldots, A_h) \subseteq [ -C n, Cn ].
\]
Therefore,
\[
\left| \varphi(A_1,\ldots, A_h) \right|  \leq 2C n+1.
\]
If $\mca = (A_1,\ldots, A_h)$ is a $\varphi$-Sidon system of multiplicity $g$, then 
\[
g\left| \varphi(A_1,\ldots, A_h) \right| \geq  \prod_{i=1}^h |A_i|.   
\]
If $|A_i| = q$ for all $i \in \{1,\ldots, h\}$, then 
\beq                                                            \label{BC-Sidon:qCn}     
q^h =  \prod_{i=1}^h |A_i| \leq g(2Cn+1). 
\eeq

For every positive integer $g$, let $F_{\varphi,g}(n)$ denote the largest integer $q$ such that there exists 
a $\varphi$-Sidon system $\mca = (A_1,\ldots, A_h)$ of multiplicity $g$ such that 
\[
A_i \subseteq [1,n] \qqand |A_i| = q
\]
for all $i =1,\ldots, h$.    Inequality~\eqref{BC-Sidon:qCn}      implies that 
\beq                                                            \label{BC-Sidon:FgCn}      
\limsup_{n\rightarrow \infty} \frac{F_{\varphi,g}(n)}{n^{1/h}} \leq \left(2gC\right)^{1/h}  < \infty.
\eeq


\section{Constructing large $\varphi$-Sidon systems} 

\bl                                                     \label{BC-Sidon:lemma:uh-1} 
Let $p$ be a prime number and let $h \geq 2$.  There exists an integer $u_p$ such that 
\beq                    \label{BC-Sidon:uh-1} 
\gcd(u_p,p) =  \gcd(u_p^h -1,p) = 1
\eeq
if and only if $p-1$ does not divide $h$.
\el

\begin{proof}
Because $p$ is prime we have $\gcd(u,p) = 1 \text{ or } p$ for all integers $u$.  
If $p-1$ divides $h$, then $h = (p-1)t$ for some integer $t$.  For all integers $u \not\equiv 0 \pmod{p}$ 
we have $u^{p-1} \equiv 1 \pmod{p}$ and so 
\[
u^h \equiv \left(u^{p-1}\right)^t \equiv 1 \pmod{p}
\]
and  
\[
\gcd(u^h -1,p) = p.
\]
Thus, if $p-1$ divides $h$, then no integer satisfies~\eqref{BC-Sidon:uh-1}.

Suppose that $p-1$ does not divide $h$. 
For every primitive root $u_p$ modulo $p$ we have
 $u_p \not\equiv 0 \pmod{p}$  and  $u_p^h \not\equiv 1 \pmod{p}$ 
and so $u_p$ satisfies~\eqref{BC-Sidon:uh-1}.
This completes the proof. 
\end{proof}

\bl                                                \label{BC-Sidon:lemma:qh-1} 
Let $h \geq 2$ and let $\mcp(h)$ be the finite set of primes $p$ such that $p-1$ divides $h$.  
Let $\varphi = c_1x_1+\cdots + c_hx_h$ be a linear form whose coefficients $c_1,\ldots, c_h$ 
are nonzero integers such that $\gcd(c_i, p)=1$ for all $i \in \{1,\ldots, h\}$ and $p \in \mcp(h)$.  
There exist positive integers $u$ and $Q$ with $\gcd(u,Q) = 1$ such that every prime 
$q$ in the infinite arithmetic progression $u \pmod{Q}$ satisfies  
\[
\gcd(q^h -1,c_i) = 1
\]
for all $i \in \{1,\ldots, h\}$. 
\el

\begin{proof}
Let $\mcp(c_1,\ldots, c_h)$ be the finite set of primes $p$ 
such that $p$ divides $c_i$ for some $i \in \{1,\ldots, h\}$ and let 
\[
Q = \prod_{p \in \mcp(c_1,\ldots, c_h)} p.  
\]   
If $p \in \mcp(c_1,\ldots, c_h)$, then $p \notin \mcp(h)$ and so $p-1$ does not divide $h$,  
and so, by Lemma~\ref{BC-Sidon:lemma:uh-1}, there is an integer $u_p$ that satisfies~\eqref{BC-Sidon:uh-1}.
By the Chinese remainder theorem, there is an integer $u$ such that 
\[
u \equiv u_p \pmod{p}
\]
for all $p \in \mcp(c_1,\ldots, c_h)$.  
We have $\gcd(u,p) = \gcd(u_p,p) = 1$ for all $p \in \mcp(c_1,\ldots, c_h)$ and so  $\gcd(u,Q)=1$.  
By Dirichlet's theorem, there are infinitely many primes $q$ 
such that 
\[
q \equiv u \pmod{Q}. 
\]
For all $p \in \mcp(c_1,\ldots, c_h)$ we have
\[
\gcd(q^h -1,p) = \gcd(u^h -1,p) = \gcd(u_p^h -1,p) =  1 
\]
and so 
\[
\gcd(q^h -1,c_i) = 1 
\]
for all $i \in \{1,\ldots, h\}$. 
This completes the proof.  
\end{proof}

\bt                                                        \label{BC-Sidon:theorem:BCN-1}     
Consider the linear form 
\[
\varphi = c_1 x_1 + \cdots + c_h x_h
\]
where $h \geq 2$ and $c_1,\ldots, c_h$ are nonzero integers.  
Let $q$ be a prime such that 
\[
\gcd(q^h-1,c_i) = 1
\]
for all $i \in \{1,\ldots, h\}$.
There is an $h$-tuple $\mca = (A_1,\ldots, A_h)$ of sets of integers with 
\[
A_i \subseteq \{1,\ldots, q^h -2\}  \qqand   |A_i| = q 
\]
for all $i \in \{1,\ldots, h\}$ such that $\mca$ is a $\varphi$-Sidon system of multiplicity at most $h!$.  
Equivalently, 
\beq                                                            \label{BC-Sidon:BCN}     
F_{\varphi,h!} (q^h -2) \geq q. 
\eeq
\et

\begin{proof}
Let $\Fq$ be the finite field with $q$ elements and let $\F_{q^h}$ be an extension field of \Fq\ 
of degree $h$.  The multiplicative group $\F_{q^h}^{\times}$ is  cyclic 
of order $q^h - 1$.  Let $\theta$ be a generator of $\F_{q^h}^{\times}$.  
For all $i \in \{1,\ldots, h\}$ we have $\gcd(q^h-1,c_i) = 1$ 
and so $\theta^{c_i}$ is also a generator of the cyclic group $\F_{q^h}^{\times}$.  

Let $\Fq = \{\lambda_1,\ldots, \lambda_q \}$.  
The inequality $h \geq 2$ implies that $\theta \notin \F_q$ 
and so, for all $j \in \{1,\ldots, q\}$, 
we have  $\theta \neq \lambda_j$.  Equivalently, 
$\theta - \lambda_j \in \F_{q^h}^{\times}$.
For all $i \in \{1,\ldots, h\}$, the element $\theta^{c_i}$ generates $\F_{q^h}^{\times}$.
It follows that there is a unique integer 
\[
a_{i,j} \in \{0,1,\ldots, q^h -2\}
\] 
such that 
\[
\theta^{c_i a_{i,j}} = \theta - \lambda_j.
\]
If $a_{i,j} = 0$, then $\theta = \lambda_j + 1 \in \Fq$, which is absurd.  Therefore, 
\[
a_{i,j} \in \{1,\ldots, q^h -2\}
\] 
for all $i \in \{1,\ldots,h\}$ and $j \in \{1,\ldots, q \}$.
Moreover, if $1 \leq j < k \leq q$, then $\lambda_j \neq \lambda_k$ and so $a_{i,j} \neq a_{i,k}$ .  
The set 
\[
A_i = \left\{ a_{i,j}: j = 1,\ldots, q \right\}
\]   
is a subset of $\{1,\ldots, q^h -2\}$ of cardinality $q$ for all $i =1,\ldots,h$.  

Let $\mca = (A_1,\ldots, A_h)$.
For all $(a_{1,j_1},\ldots, a_{h,j_h})  \in A_1 \times \cdots \times A_h$ we have   
\begin{align*}
\theta^{\varphi \left( a_{1,j_1},\ldots, a_{h,j_h} \right) } & =  \theta^{\sum_{i=1}^h c_ia_{i,j_i}  }
 =  \prod_{i=1}^h  \theta^{c_ia_{i,j_i}  }  =  \prod_{i=1}^h (\theta - \lambda_{j_i}) = f(\theta)
\end{align*} 
where 
\[
f(t) =  \prod_{i=1}^h (t- \lambda_{j_i}) = t^h  - \left(\sum_{i=1}^h \lambda_{j_i} \right) t^{h-1} + \cdots 
+ (-1)^h \prod_{i=1}^h\lambda_{j_i} 
\]
is a monic polynomial of degree $h$ with coefficients in \Fq. 

Similarly, for $(a_{1,j'_1},\ldots, a_{h,j'_h})  \in A_1 \times \cdots \times A_h$ we have  
\begin{align*}
\theta^{\varphi\left( a_{1,j'_1},\ldots, a_{h,j'_h} \right) } 
 &  =  \prod_{i=1}^h (\theta - \lambda_{j'_i}) = g(\theta)
\end{align*} 
where 
\[
g(t) = \prod_{i=1}^h (t- \lambda_{j'_i}) = t^h - \left(\sum_{i=1}^h \lambda_{j'_i} \right) t^{h-1} + \cdots 
+ (-1)^h \prod_{i=1}^h\lambda_{j'_i} 
\]
is also a monic polynomial of degree $h$ with coefficients in \Fq. 

The relation  $\varphi\left( a_{1,j_1},\ldots, a_{h,j_h} \right) = \varphi\left( a_{1,j'_1},\ldots, a_{h,j'_h} \right)$ 
implies 
\[
f(\theta) = \theta^{\varphi \left( a_{1,j_1},\ldots, a_{h,j_h} \right) }  
=  \theta^{\varphi\left( a_{1,j'_1},\ldots, a_{h,j'_h} \right) }  = g(\theta) 
\]
and so $\theta$ is a root of the polynomial $f(t)-g(t)$.  
If $f(t) \neq g(t)$, then  $f(t) - g(t)$ is a nonzero polynomial of degree at most $h-1$.  
This is impossible because the minimal polynomial of $\theta$ has degree $h$.  
Therefore, 
\[
 \prod_{i=1}^h (t- \lambda_{j_i}) = f(t) = g(t)  = \prod_{i=1}^h (t- \lambda_{j'_i})
 \]
and so  $(\lambda_{j'_1},\ldots, \lambda_{j'_h})$ is a permutation of 
$(\lambda_{j_1},\ldots, \lambda_{j_h})$ and  $(j'_1,\ldots, j'_h)$ is a permutation of 
$(j_1,\ldots, j_h)$.   There are at most $h!$ such permutations.   
It follows that for every integer $w$ there are at most $h!$ elements 
$(a_1,\ldots, a_h) \in A_1\times \cdots \times A_h$ 
such that $\varphi(a_1,\ldots, a_h) = w$.  
Therefore, $R_{\mca,\varphi}(w) \leq h!$ and $\mca = (A_1,\ldots, A_h)$ 
is a $\varphi$-Sidon system of multiplicity at most $h!$. 
This completes the proof. 
\end{proof}

\bt                                                        \label{BC-Sidon:theorem:BCN}     
Let $h \geq 2$ and let $\mcp(h)$ be the finite set of primes $p$ such that $p-1$ divides $h$.  
Let $\varphi = c_1x_1+\cdots + c_hx_h$ be a linear form whose coefficients $c_1,\ldots, c_h$ 
are nonzero integers such that for all $i \in \{1,\ldots, h\}$ and $p \in \mcp(h)$.  Then  
\[
\liminf_{n\rightarrow\infty} \frac{F_{\varphi,h!}(n)}{n^{1/h}}
\geq 1.
\]
\et

\begin{proof}
There is an analog of Hoheisel's theorem for sufficiently large primes in arithmetic progressions.  
Let $u$ and $Q$ be relatively prime positive integers.  
Baker, Harman, and Pintz~\cite{bake-harm-pint97} proved that there 
is a real number  $\alpha = \alpha(u,Q)$ with $0 < \alpha < 1$ such that if $p$ and $p'$ 
are sufficiently large consecutive primes in the arithmetic progression 
$u \pmod{Q}$, then $p' - p < p^{\alpha}$.

By Lemma~\ref{BC-Sidon:lemma:qh-1},there exist positive integers $u$ and $Q$ with $\gcd(u,Q) = 1$ 
such that every prime 
$q$ in the infinite arithmetic progression $u \pmod{Q}$ satisfies  
\[
\gcd(q^h -1,c_i) = 1
\]
for all $i \in \{1,\ldots, h\}$. 

For every sufficiently large integer $n$, let $p$ be the largest prime such that $p\equiv u \pmod{Q}$ and 
$p \leq n^{1/h}$.  Let $p'$ be the smallest prime such that $p\equiv u \pmod{Q}$ and $p' > n^{1/h}$.  
Then $p$ and $p'$ are consecutive  primes in the arithmetic progression  $ u \pmod{Q}$, and so 
\[
p \leq n^{1/h} < p' \leq p+p^{\alpha} \leq p+  n^{\alpha/h}. 
\]
Applying inequality~\eqref{BC-Sidon:BCN} from Theorem~\ref{BC-Sidon:theorem:BCN-1}   
with $q = p$, we obtain 
\[
F_{\varphi,h!}(n)  \geq F_{\varphi,h!}(p^h)  \geq F_{\varphi,h!}(p^h-2) \geq p 
\geq n^{1/h} - n^{\alpha /h}.
\]
Therefore,
\[
\liminf_{n\rightarrow \infty} \frac{F_{\varphi,h!}(n )}{n^{1/h}} 
\geq \liminf_{n\rightarrow \infty} \frac{ n^{1/h} - n^{\alpha /h}}{n^{1/h}} 
= \liminf_{n\rightarrow \infty} \left( 1 - \frac{ 1}{ n^{(1-\alpha) /h}} \right) = 1.
\]
This completes the proof. 
\end{proof}


\section{Open problems}
Let  $h \geq 2$ and let $\varphi = \sum_{i=1}^h c_i x_i$ be a  linear form 
with nonzero integer coefficients.  
\benum
\item[(i)]
Is it true that 
\[
\liminf_{n\rightarrow\infty} \frac{F_{\varphi,h!}(n)}{n^{1/h}} > 0
\]
with no condition on the primes that divide the coefficients $c_i$? 
\item[(ii)]
Is it true that there is a finite set $\mcp(h)$ of prime numbers 
such that if none  of  the coefficients of $\varphi$ is divisible by a prime 
in $\mcp(h)$, then 
\[
\liminf_{N\rightarrow \infty} \frac{F_{\varphi}(N)}{N^{1/h}}  > 0?
\]
\item[(iii)]
Is it true that 
\[
\liminf_{N\rightarrow \infty} \frac{F_{\varphi}(N)}{N^{1/h}}  > 0 
\]
with no condition on the primes that divide the coefficients $c_i$? 
This would be the analog of the Bose-Chowla theorem for classical Sidon sets.  
\eenum

\def\cprime{$'$} \def\cprime{$'$}
\providecommand{\bysame}{\leavevmode\hbox to3em{\hrulefill}\thinspace}
\providecommand{\MR}{\relax\ifhmode\unskip\space\fi MR }
\providecommand{\MRhref}[2]{%
  \href{http://www.ams.org/mathscinet-getitem?mr=#1}{#2}
}
\providecommand{\href}[2]{#2}


\begin{thebibliography}{1}

\bibitem{bake-harm-pint97}
R.~C. Baker, G.~Harman, and J.~Pintz, \emph{The exceptional set for
  {G}oldbach's problem in short intervals}, Sieve methods, exponential sums,
  and their applications in number theory ({C}ardiff, 1995), London Math. Soc.
  Lecture Note Ser., vol. 237, Cambridge Univ. Press, Cambridge, 1997,
  pp.~1--54.

\bibitem{bose-chow62}
R.~C. Bose and S.~Chowla, \emph{Theorems in the additive theory of numbers},
  Comment. Math. Helv. \textbf{37} (1962/63), 141--147.

\bibitem{halb-roth66}
H.~Halberstam and K.~F. Roth, \emph{{Sequences, Vol. 1}}, Oxford University
  Press, Oxford, 1966, Reprinted by Springer-Verlag, Heidelberg, in 1983.

\bibitem{heat88b}
D.R. Heath-Brown, \emph{The number of primes in a short interval}, {J. reine
  angew. Math} \textbf{389} (1988), 22--63.

\bibitem{hohe30}
G.~Hoheisel, \emph{{Primzahlprobleme in der Analysis}}, {Sitz. Preuss. Akad.
  Wiss.} \textbf{2} (1930), 1--13.

\bibitem{obry04}
K.~{O'Bryant}, \emph{A complete annotated bibliography of work related to
  {Sidon} sequences}, Electronic J. Combinatorics (2004), Dynamic Surveys DS
  11.

\end{thebibliography}
\end{document}